\documentclass[11pt,a4paper]{article}
\usepackage{bbding}
\usepackage{amsfonts}
\usepackage{epsf,epsfig,amsfonts,amsgen,amsmath,amstext,amsbsy,amsopn,amsthm,lineno}
\usepackage{color}
\newtheorem{theorem}{Theorem}

\newtheorem{false statement}{False Statement}

\theoremstyle{definition}

\newtheorem{problem}{Problem}
\newtheorem{case}{Case}
\newtheorem{subcase}{Case}[case]
\newtheorem{subsubcase}{Case}[subcase]
\newtheorem{subsubsubcase}{Case}[subsubcase]

\begin{document}

\title{\bf\Large On heavy paths in 2-connected weighted graphs\thanks{Supported  by NSFC
(No.~10871158).}}

\date{}

\author{Binlong Li ~and Shenggui
Zhang\thanks{Corresponding author.
E-mail address: sgzhang@nwpu.edu.cn}\\[2mm]
\small Department of Applied Mathematics,
\small Northwestern Polytechnical University,\\
\small Xi'an, Shaanxi 710072, P.R.~China\\}
\maketitle

\begin{abstract}
A weighted graph is a graph in which every edge is assigned a
non-negative real number. In a weighted graph, the weight of a path
is the sum of the weights of its edges, and the weighed degree of a
vertex is the sum of the weights of the edges incident with it. In
this paper we give three weighted degree conditions for the
existence of heavy or Hamilton paths with one or two given
end-vertices in 2-connected weighted graphs.

\medskip
\noindent {\bf Keywords:} Weighted graph; Heavy path; Weighed degree
\smallskip
\end{abstract}

\section{Introduction}

We use Bondy and Murty \cite{Bondy_Murty} for terminology and
notation not defined here and consider finite simple graphs only.

A {\it weighted graph} is a graph in which every edge $e$ is
assigned a non-negative real number $w(e)$, called the {\it weight}
of $e$. Let $G=(V, E)$ be a weighted graph. For a subgraph $H$ of
$G$, $V(H)$ and $E(H)$ denote the sets of vertices and edges of $H$,
respectively. The {\it weight} of $H$ is defined by
$$
w(H)=\sum_{e\in E(H)}w(e).
$$
For a vertex $v\in V$, $N_H(v)$ denotes the set, and $d_H(v)$ the
number, of vertices in $H$ that are adjacent to $v$. We define the
{\it weighted degree} of $v$ in $H$ by
$$
d_H^w(v)=\sum_{x\in N_H(v)}w(vx).
$$
When no confusion occurs, we will denote $N_G(v)$, $d_G(v)$ and
$d_G^w(v)$ by $N(v)$, $d(v)$ and $d^w(v)$, respectively.

An unweighted graph can be regarded as a weighted graph in which
each edge $e$ is assigned weight $w(e)=1$. Thus, in an unweighted
graph, $d^w(v)=d(v)$ for every vertex $v$, and the weight of a
subgraph is simply the number of its edges.

About twenty years ago, Bondy and Fan \cite{Bondy_Fan} began the
study of the existence of heavy paths and cycles in weighted graphs.
They showed two results with Dirac-type weighted degree condition.
In the following by an {\it $x$-path} we mean a path whose initial
vertex is $x$; and by an {\it $(x,y)$-path} we mean one whose
end-vertices are $x$ and $y$.

\begin{theorem}[Bondy and Fan \cite{Bondy_Fan}]
Let $G$ be a 2-connected weighted graph, $x$ and $y$ be two vertices
of $G$, and $d$ be a real number. If $d^w(v)\geq d$ for every vertex
$v$ in $V(G)\backslash\{x,y\}$, then $G$ contains an $(x,y)$-path of
weight at least $d$.
\end{theorem}

\begin{theorem}[Bondy and Fan \cite{Bondy_Fan}]
Let $G$ be a 2-connected weighted graph and $d$ be a real number. If
$d^w(v)\geq d$ for every vertex $v$ in $V(G)$, then either $G$
contains a cycle of weight at least $2d$ or every heaviest cycle in
$G$ is a Hamilton cycle.
\end{theorem}

Corresponding Ore-type conditions for the existence of heave paths
and cycles in weighted graphs were also obtained.

\begin{theorem}[Enomoto, Fujisawa and Ota \cite{Enomoto_Fujisawa_Ota}]
Let $G$ be a 2-connected weighted graph, $x$ and $y$ be two vertices
of $G$, and $d$ be a real number. If $d^w(v_1)+d^w(v_2)\geq 2d$ for
every pair of nonadjacent vertices $v_1$ and $v_2$ in
$V(G)\backslash\{x,y\}$, then $G$ contains either an $(x,y)$-path of
weight at least $d$ or a Hamilton $(x,y)$-path.
\end{theorem}

\begin{theorem}[Bondy~{\em et al.}~\cite{BBHV}]
Let $G$ be a 2-connected weighted graph and $d$ be a real number. If
$d^w(v_1)+d^w(v_2)\geq 2d$ for every pair of nonadjacent vertices
$v_1$ and $v_2$ in $V(G)$, then $G$ contains either a cycle of
weight at least $2d$ or a Hamilton cycle.
\end{theorem}

It can be seen that Theorems 2 and 4 generalize the classical
results of Dirac and Ore on the existence of long cycles in
unweighted graphs, respectively. Another well-known result on the
existence of long cycles in unweighted graphs is the so-called
Fan-type condition \cite{Fan}. By constructing some examples,
Zhang~{\em et al.}~\cite{ZBLW} found that the Fan's condition cannot
be generalized to weighted graphs directly. As one can expect, the
corresponding condition for heavy paths in weighted graphs is also
not valid.

\begin{false statement}
Let $G$ be a 2-connected weighted graph, $x$ and $y$ be two vertices
of $G$, and $d$ be a real number. If $\max\{d^w(v_1),d^w(v_2)\}\geq
d$ for every pair of vertices $v_1$ and $v_2$ in
$V(G)\backslash\{x,y\}$ with distance 2, then $G$ contains either an
$(x,y)$-path of weight at least $d$ or a Hamilton $(x,y)$-path.
\end{false statement}

This can be shown by the graph in Fig.~1. We assign weight 2 to the
edge $u_2u_3$, and weight 1 to all the remaining edges. Then, the
resulting weighted graph is 2-connected and
$\max\{d^w(v_1),d^w(v_2)\}\geq 5$ for every pair of vertices $v_1$
and $v_2$ in $V(G)\backslash\{x,y\}$ with distance 2. Whereas, the
graph contains neither an $(x,y)$-path of weight at least 5 nor a
Hamilton $(x,y)$-path.

\begin{center}

\begin{picture}(80,60)

\put(40,0){\circle*{4}} \put(0,20){\circle*{4}}
\put(80,20){\circle*{4}} \put(0,40){\circle*{4}}
\put(80,40){\circle*{4}} \put(40,60){\circle*{4}}

\put(40,0){\line(-2,1){40}} \put(40,0){\line(2,1){40}}
\put(0,20){\line(1,0){80}} \put(0,20){\line(0,1){20}}
\put(0,20){\line(4,1){80}} \put(80,20){\line(-4,1){80}}
\put(80,20){\line(0,1){20}} \put(0,40){\line(2,1){40}}
\put(80,40){\line(-2,1){40}}

\put(-7.5,37.5){$x$} \put(82.5,37.5){$y$} \put(37.5,62.5){$u_1$}
\put(-10,17.5){$u_2$} \put(82.5,17.5){$u_3$} \put(37.5,-7.5){$u_4$}

\end{picture}\\

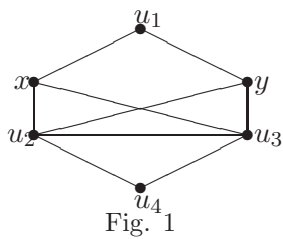
\small{Fig. 1}

\end{center}

Here we strengthen the condition of False Statement 1 and get the
following result:

\begin{theorem}
Let $G$ be a 2-connected weighted graph, $x$ and $y$ be two vertices
of $G$, and $d$ be a real number. If $\max\{d^w(v_1),d^w(v_2)\}\geq
d$ for every pair of nonadjacent vertices $v_1$ and $v_2$ in
$V(G)\backslash\{x,y\}$, then $G$ contains either an $(x,y)$-path of
weight at least $d$ or a Hamilton $(x,y)$-path.
\end{theorem}

In fact, a corresponding result on the existence of heavy cycles in
weighted graphs was obtained by Fujisawa \cite{Fujisawa}.

\begin{theorem}[Fujisawa \cite{Fujisawa}]
Let $G$ be a 2-connected weighted graph and $d$ be a real number. If
$\max\{d^w(v_1),d^w(v_2)\}\geq d$ for every pair of nonadjacent
vertices $v_1$ and $v_2$ in $V(G)$, then $G$ contains either a cycle
of weight at least $2d$ or a Hamilton cycle.
\end{theorem}

If we weaken the weighted degree sum condition of two nonadjacent
vertices in Theorem 3 by that of three pairwise nonadjacent
vertices, can we still have the same result? The answer is also
negative.

\begin{false statement}
Let $G$ be a 2-connected weighted graph, $x$ and $y$ be two vertices
of $G$, and $d$ be a real number. If $d^w(v_1)+d^w(v_2)+d^w(v_3)\geq
3d$ for every three pairwise nonadjacent vertices $v_1$, $v_2$ and
$v_3$ in $V(G)\backslash\{x,y\}$, then $G$ contains either an
$(x,y)$-path of weight at least $d$ or a Hamilton $(x,y)$-path.
\end{false statement}

This can be shown as follows. Let $H_1$ and $H_2$ be two complete
graphs on at least three vertices with two vertices $x$ and $y$ in
common. Let $G=H_1\cup H_2$ and assign weight 0 to every edge of
$G$. Then $G$ is 2-connected and there exist no three pairwise
nonadjacent vertices in $G$. So $G$ satisfies the condition of False
Statement 2 for any $d>0$, but contains neither an $(x,y)$-path of
weight at least $d$ nor a Hamilton $(x,y)$-path.

When specifying only one end-vertex of the path, we have the
following result:

\begin{theorem}
Let $G$ be a 2-connected weighted graph, $x$ be a vertex of $G$, and
$d$ be a real number. If $d^w(v_1)+d^w(v_2)+d^w(v_3)\geq 3d$ for
every three pairwise nonadjacent vertices $v_1$, $v_2$ and $v_3$ in
$V(G)\backslash\{x\}$, then $G$ contains either an $x$-path of
weight at least $d$ or a Hamilton $x$-path.
\end{theorem}

Instead of Theorem 7, we will prove the following stronger result.

\begin{theorem}
Let $G$ be a 2-connected weighted graph, $x$ be a vertex of $G$, and
$d$ be a real number. If $\max\{d^w(v_1),d^w(v_2),d^w(v_3)\}\geq d$
for every three pairwise nonadjacent vertices $v_1$, $v_2$ and $v_3$
in $V(G)\backslash\{x\}$, then $G$ contains either an $x$-path of
weight at least $d$ or a Hamilton $x$-path.
\end{theorem}

Further, under the condition of False Statement 2, we can prove the
following result:

\begin{theorem}
Let $G$ be a 2-connected weighted graph, $x$ and $y$ be two vertices
of $G$, and $d$ a real number. If $d^w(v_1)+d^w(v_2)+d^w(v_3)\geq
3d$ for every three pairwise nonadjacent vertices $v_1$, $v_2$ and
$v_3$ in $V(G)\backslash\{x,y\}$, then at least one of the following
holds:\\
$(1)$ $G$ contains an $(x,y)$-path $P$ with $w(P)\geq d$;\\
$(2)$ $G$ contains an $x$-path $P_1$ and a $y$-path $P_2$
which are disjoint with  $w(P_1)+w(P_2)\geq d$;\\
$(3)$ $G$ contains an $x$-path $P_1$ and a $y$-path $P_2$ which are
disjoint with  $V(G)=V(P_1)\cup V(P_2)$.
\end{theorem}

Instead of Theorem 9, we will prove the following stronger result.

\begin{theorem}
Let $G$ be a 2-connected weighted graph, $x$ and $y$ be two vertices
of $G$, and $d$ be a real number. If
$\max\{d^w(v_1),d^w(v_2),d^w(v_3)\}\geq d$ for every three pairwise
nonadjacent vertices $v_1$, $v_2$ and $v_3$ in
$V(G)\backslash\{x,y\}$, then at least one of the following
holds:\\
$(1)$ $G$ contains an $(x,y)$-path $P$ with $w(P)\geq d$;\\
$(2)$ $G$ contains an $x$-path $P_1$ and a $y$-path $P_2$ which are
disjoint with $w(P_1)+w(P_2)\geq d$;\\
$(3)$ $G$ contains an $x$-path $P_1$ and a $y$-path $P_2$ which are
disjoint with  $V(G)=V(P_1)\cup V(P_2)$.
\end{theorem}

We give the proofs of Theorems 5, 8 and 10 in the following
sections.

\section{Proof of Theorem 5}

If $d=0$, then the assertion is obvious. Hence we may assume that
$d>0$. Let $|V(G)|=n$. We use induction on $n$.

If $n=3$, then the result is trivially true. Suppose now that $n\geq
4$ and the theorem is true for all graphs on fewer than $n$
vertices. Let $H=G-x$.

\begin{case}
$H$ is 2-connected.
\end{case}

Since $G$ is 2-connected, we have $d(x)\geq 2$. Choose a vertex
$x'\in N(x)\backslash\{y\}$ such that $w(xx')=\max\{w(xv): v\in
N(x)\backslash\{y\}\}$. Then for every $v\in V(H)\backslash\{y\}$,
we have $d^w_H(v)\geq d^w_G(v)-w(xx')$. Hence,
$$
\max\{d^w_H(v_1),d^w_H(v_2)\}\geq d-w(xx')
$$
for every pair of nonadjacent vertices $v_1$ and $v_2$ in
$V(H)\backslash\{x',y\}$. By the induction hypothesis, there is an
$(x',y)$-path $P'$ in $H$ such that either $w(P')\geq d-w(xx')$ or
$P'$ is a Hamilton path of $H$. Thus the path $P=xx'P'$ is a
required path.

\begin{case}
$H$ is separable.
\end{case}

Let $z$ be a cut vertex of $H$, $H_1$ be a component of $H-z$, and
$H_2=H-z-H_1$. For $i=1,2$, let $G_i=G[H_i\cup\{x,z\}]$. If
$xz\notin E(G_i)$, we add the edge $xz$ of weight 0 to $G_i$. Now
the two resulting graphs $G_1$ and $G_2$ are both 2-connected.

\begin{subcase}
$y=z$.
\end{subcase}

Let $d^w(v_0)=\min\{d^w(v): v\in V(H_1)\cup V(H_2)\}$. Without loss
of generally, we can assume that $v_0\in V(H_1)$. Since for all
$v\in V(H_2)$, $vv_0\notin E(G)$, we have
$d^w(v)=\max\{d^w(v),d^w(v_0)\}\geq d$. Then
$d^w_{G_2}(v)=d^w(v)\geq d$ for all $v\in V(G_2)\backslash\{x,y\}$.
By Theorem 1, there exists an $(x,y)$-path $P$ in $G_2$ such that
$w_{G_2}(P)\geq d$. It is obvious that $P$ is not the added edge
$xy$ itself, so we can take the path $P$ as a required path.

\begin{subcase}
$y\neq z$.
\end{subcase}

Without loss of generally, we can assume that $y\in V(H_1)$. It is
easy to know that $\max\{d^w_{G_2}(v_1),d^w_{G_2}(v_2)\}\geq d$ for
every pair of nonadjacent vertices in $V(G_2)\backslash\{x,z\}$. By
the induction hypothesis, there exists an $(x,z)$-path $P_2$ in
$G_2$ such that either $w_{G_2}(P_2)\geq d$ or $P_2$ is a Hamilton
path in $G_2$. It is obvious that $P_2$ is not the added edge $xy$
itself.

Suppose that $w_{G_2}(P_2)\geq d$. Let $P_1$ be a $(z,y)$-path not
passing through $x$ in $G_1$. Then the path $P=P_2P_1$ is a required
path.

Suppose now that every $(x,z)$-path in $G_2$ has weight less than
$d$. Then $P_2$ is a Hamilton path in $G_2$, and it follows from
Theorem 1 that there exists at least one vertex in $H_2$ with
weighted degree less than $d$. Thus, every vertex in
$V(H_1)\backslash\{y\}$ has weighted degree at least $d$.

\begin{subsubcase}
$H_2$ has at least two vertices.
\end{subsubcase}

Let $G'_1$ be the weighted graph such that
$V(G'_1)=V(G_1)\cup\{x'\}$, where $x'\notin V(G)$;
$E(G'_1)=E(G_1)\cup\{xx',x'z\}$; and
$$
w_{G'_1}(e)=\left\{
\begin{array}{ll}
  w_{G_1}(e),   & \mbox{if\ } e\in E(G_1);\\
  0,            & \mbox{otherwise}.
\end{array}
\right.
$$
It is easy to know that $G'_1$ is 2-connected. Besides, we can see
that each vertex in $V(G'_1)\backslash\{x,y\}$ other than $x'$ and
$z$ has weighted degree at least $d$. Thus,
$$
\max\{d^w_{G'_1}(v_1),d^w_{G'_1}(v_2)\}\geq d
$$
for every pair of nonadjacent vertices $v_1$ and $v_2$ in
$V(G'_1)\backslash\{x,y\}$. By the induction hypothesis, there
exists an $(x,y)$-path $P'$ in $G'_1$ such that either
$w_{G'_1}(P')\geq d$ or $P'$ is a Hamilton path in $G'_1$.

Suppose that $w_{G'_1}(P')\geq d$. If $P'$ contains neither the
added edge $xz$ nor the subpath $xx'z$, then it is a required path.
Otherwise, we can use any $(x,z)$-path in $G_2$ to replace the added
edge $xz$ or the subpath $xx'z$, and then obtain an $(x,y)$-path in
$G$ with weight at least $w_{G'_1}(P')\geq d$.

Suppose that $P'$ is a Hamilton path in $G'_1$. Because
$d_{G'_1}(x')=2$, $P'$ contains the subpath $xx'z$. Thus, we can use
the path $P_2$ to replace the subpath $xx'z$, and then obtain a
Hamilton $(x,y)$-path in $G$. \vspace{2mm}

\begin{subsubcase}
$H_2$ has only one vertex.
\end{subsubcase}

Let $x'$ be the vertex in $H_2$. Then we have $x'\in
N(x)\backslash\{y\}$, $d(x')=2$ and $d^w(x')<d$. By a similar proof,
we can obtain that there exists a vertex $y'$ in
$N(y)\backslash\{x\}$ such that $d(y')=2$ and $d^w(y')<d$. Since for
all vertex $v\in V(H_1)\backslash\{y\}$, $d^w(v)\geq d$, we have
$y'=z$ and $d(z)=2$.

If $H_1$ has only one vertex $y$, then $P=xx'zy$ is a Hamilton
$(x,y)$-path in $G$. Thus, we can assume that $H_1$ has at least one
vertex other than $y$. Let $G'_1=G[H_1\cup\{x\}]$. If $xy\notin
E(G'_1)$, we add the edge $xy$ of weight zero to $G'_1$, and then
get that $G'_1$ is 2-connected and $w_{G'_1}(v)=w_G(v)\geq d$ for
every vertex $v\in V(G'_1)\backslash\{x,y\}$. By Theorem 1, there
exists an $(x,y)$-path $P$ in $G'_1$ such that $w_{G'_1}(P)\geq d$.
It is obvious that $P$ is not the added edge $xy$ itself. So we can
take the path $P$ as a required path.

The proof is complete.\hfill$\Box$

\section{Proof of Theorem 8}

If $d=0$, then the assertion is obvious. Hence we may assume that
$d>0$. Let $|V(G)|=n$. We use induction on $n$.

If $n=3$, then the result is trivially true. Suppose now that $n\geq
4$ and the theorem is true for all graphs on fewer than $n$
vertices. Let $H=G-x$.

\setcounter{case}{0}
\begin{case} $H$ is
2-connected.
\end{case}

Since $G$ is 2-connected, we have $d(x)\geq 2$. Choose $x'\in N(x)$
such that $w(xx')=\max\{w(xv): v\in N(x)\}$. Then for every $v\in
V(H)$, we have $d^w_H(v)\geq d^w_G(v)-w(xx')$. Hence,
$$
\max\{d^w_H(v_1),d^w_H(v_2),d^w_H(v_3)\}\geq d-w(xx')
$$
for every three nonadjacent vertices $v_1$, $v_2$ and $v_3$ in
$V(H)\backslash\{x'\}$. By the induction hypothesis, there is an
$x'$-path $P'$ in $H$ such that either $w(P')\geq d-w(xx')$ or $P'$
is a Hamilton path of $H$. Then the path $P=xx'P'$ is a required
path.

\begin{case}
$H$ is separable.
\end{case}

Let $y$ be a cut vertex of $H$, $H_1$ be a component of $H-y$, and
$H_2=H-y-H_1$. For $i=1,2$, let $G_i=G[H_i\cup\{x,y\}]$. If
$xy\notin E(G_i)$, we add the edge $xy$ of weight 0 to $G_i$. Now
the two resulting graphs $G_1$ and $G_2$ are both 2-connected.

If for some $i\in\{1,2\}$, $d^w_{G_i}(v)=d^w(v)\geq d$ for all $v$
in $V(G_i)\backslash\{x,y\}$, then by Theorem 1, there is an
$(x,y)$-path $P$ such that $w(P)\geq d$. It is obvious that $P$ is
not the added edge $xy$ itself, so we can take the path $P$ as a
required path.

Otherwise, for $i\in\{1,2\}$, there exists a vertex $v$ in
$V(G_i)\backslash\{x,y\}$ such that $d^w_{G_i}(v)=d^w(v)<d$. Then
for $i\in\{1,2\}$, $\max\{d^w_{G_i}(v_1),d^w_{G_i}(v_2)\}\geq d$ for
every pair of nonadjacent vertices $v_1$ and $v_2$ in
$V(G_i)\backslash\{x,y\}$. By Theorem 5, there exists an
$(x,y)$-path $P_i$ such that either $w_{G_i}(P_i)\geq d$ or $P_i$ is
a Hamilton path of $G_i$. It is obvious that $P_i$ is not the added
edge $xy$ itself.

If for some $i\in\{1,2\}$, $w_{G_i}(P_i)\geq d$, then we can take
the path $P_i$ as a required path. Otherwise, $P_i$ is a Hamilton
$(x,y)$-path of $G_i$, and $C=P_1P_2$ is a Hamilton cycle of $G$.
Then $G$ contains a Hamilton $x$-path.

The proof is complete.\hfill$\Box$

\section{Proof of Theorem 10}

If $d=0$, then the assertion is obvious. Hence we may assume that
$d>0$. Let $|V(G)|=n$. We use induction on $n$.

If $n=3$, then the result is trivially true. Suppose now that $n\geq
4$ and the theorem is true for all graphs on fewer than $n$
vertices. Let $H=G-x$.

\setcounter{case}{0}
\begin{case} $H$ is
2-connected.
\end{case}

Since $G$ is 2-connected, we have $d(x)\geq 2$. Choose a vertex
$x'\in N(x)\backslash\{y\}$ such that $w(xx')=\max\{w(xv): v\in
N(x)\backslash\{y\}\}$. Then for every $v\in V(H)\backslash\{y\}$,
we have $d^w_H(v)\geq d^w_G(v)-w(xx')$. Hence,
$$
\max\{d^w_H(v_1),d^w_H(v_2),d^w_H(v_3)\}\geq d-w(xx')
$$
for every three pairwise nonadjacent vertices $v_1$, $v_2$ and $v_3$
in $V(H)\backslash\{x',y\}$. By the induction hypothesis, $H$
satisfies the conclusion of the theorem. If $H$ contains an
$(x',y)$-path $P'$ such that $w(P')\geq d-w(xx')$, then $P=xx'P'$ is
an $(x,y)$-path in $G$ with weight at least $d$. If $G$ contains an
$x'$-path $P'_1$ and a $y$-path $P'_2$ which are disjoint and have
weight sum at least $d-w(xx')$, then the $x$-path $P_1=xx'P'_1$ and
the $y$-path $P_2=P'_2$ are disjoint and have weight sum at least
$d$. If $G$ contains an $x'$-path $P'_1$ and a $y$-path $P'_2$ which
are disjoint and contain all vertices of $H$, then the $x$-path
$P_1=xx'P'_1$ and the $y$-path $P_2=P'_2$ are disjoint and contain
all vertices of $G$.

\begin{case}
$H$ is separable.
\end{case}

Let $z$ be a cut vertex of $H$, $H_1$ be a component of $H-z$, and
$H_2=H-z-H_1$. For $i=1,2$, let $G_i=G[H_i\cup\{x,z\}]$. If
$xz\notin E(G_i)$, we add the edge $xz$ of weight 0 to $G_i$. Now
the two resulting graphs $G_1$ and $G_2$ are both 2-connected.

\begin{subcase}
$y=z$.
\end{subcase}

If for some $i\in\{1,2\}$, $d^w_{G_i}(v)=d^w(v)\geq d$ for all $v$
in $V(G_i)\backslash\{x,y\}$, then by Theorem 1, there is an
$(x,y)$-path $P$ such that $w_{G_i}(P)\geq d$. It is obvious that
$P$ is not the added edge $xy$ itself, so we can take the path $P$
as a required path.

Otherwise, for $i\in\{1,2\}$, there exists a vertex $v$ in
$V(G_i)\backslash\{x,y\}$ such that $d^w_{G_i}(v)=d^w(v)<d$. Then
for $i\in\{1,2\}$, we have
$\max\{d^w_{G_i}(v_1),d^w_{G_i}(v_2)\}\geq d$ for every pair of
nonadjacent vertices $v_1$ and $v_2$ in $V(G_i)\backslash\{x,y\}$.
By Theorem 5, there exists an $(x,y)$-path $P'_i$ such that either
$w_{G_i}(P'_i)\geq d$ or $P'_i$ is a Hamilton path of $G_i$. It is
obvious that $P'_i$ is not the added edge $xy$ itself.

If for some $i\in\{1,2\}$, $w_{G_i}(P'_i)\geq d$, we can take the
path $P'_i$ as a required path. Otherwise, $P'_i$ is a Hamilton
$(x,y)$-path of $G_i$, and $C=P'_1P'_2$ is a Hamilton cycle of $G$.
Thus $G$ contains an $x$-path and a $y$-path which are disjoint and
contain all vertices of $G$.

\begin{subcase}
$y\neq z$.
\end{subcase}

Without loss of generally, we can assume that $y\in V(H_1)$. It is
easy to know that
$\max\{d^w_{G_2}(v_1),d^w_{G_2}(v_2),d^w_{G_2}(v_3)\}\geq d$ for
every three pairwise nonadjacent vertices in
$V(G_2)\backslash\{x,z\}$. By the induction hypothesis, $G_2$
satisfies the conclusion of the theorem.

If $G_2$ contains an $(x,z)$-path $P'_2$ with weight at least $d$,
then let $P'_1$ be a $(z,y)$-path not passing through $x$ in $G_1$.
Then the path $P=P'_2P'_1$ is an $(x,y)$-path in $G$ with weight at
least $w(P'_2)\geq d$. So we assume that every $(x,z)$-path in $G_2$
has weight less than $d$. Then, by Theorem 1, there exists at least
one vertex in $V(H_2)$ with weighted degree less than $d$.

If $G_2$ contains an $x$-path $P'_{21}$ and a $z$-path $P'_{22}$
which are disjoint and have weight sum at least $d$, then let $P'_1$
be a $(y,z)$-path not passing through $x$ in $G_1$. So the $x$-path
path $P_1=P'_{21}$ and the $y$-path $P_2=P'_1P'_{22}$ are disjoint
and have weight sum at least $d$.

So now we assume that $G_2$ contains an $x$-path $P'_{21}$ and a
$z$-path $P'_{22}$ which are disjoint and contain all vertices of
$G_2$.

\begin{subsubcase}
$d^w_{G_1}(v)\geq d$ for all $v$ in $V(H_1)\backslash\{y\}$.
\end{subsubcase}

Let $G'_1$ be the weighted graph such that
$V(G'_1)=V(G_1)\cup\{x'\}$, where $x'\notin V(G)$;
$E(G'_1)=E(G_1)\cup\{xx',x'z\}$; and
$$
w_{G'_1}(e)=\left\{
\begin{array}{ll}
  w_{G_1}(e),   & \mbox{if\ } e\in E(G_1);\\
  0,            & \mbox{otherwise}.
\end{array}
\right.
$$
It is easy to know that $G'_1$ is 2-connected. Besides, we can see
that each vertex in $V(G'_1)\backslash\{x,y\}$ other than $x'$ and
$z$ has weighted degree at least $d$. Thus,
$$
\max\{d^w_{G'_1}(v_1),d^w_{G'_1}(v_2)\}\geq d
$$
for every pair of nonadjacent vertices $v_1$ and $v_2$ in
$V(G'_1)\backslash\{x,y\}$. By Theorem 5, there exists an
$(x,y)$-path $P'$ in $G'_1$ such that either $w_{G'_1}(P')\geq d$ or
$P'$ is a Hamilton path in $G'_1$.

Suppose that $w_{G'_1}(P')\geq d$. If $P'$ contains neither the
added edge $xz$ nor the subpath $xx'z$, then it is an $(x,y)$-path
in $G$ with weight at least $d$. Otherwise, we can use any
$(x,z)$-path in $G_2$ to replace the added edge $xz$ or the subpath
$xx'z$, and then obtain an $(x,y)$-path in $G$ with weight at least
$w_{G'_1}(P')\geq d$.

Suppose that $P'$ is a Hamilton path in $G'_1$. Because
$d_{G'_1}(x')=2$, $P'$ contains the subpath $xx'z$. Let
$P'_1=P'-\{xx',x'z\}$. So the $x$-path $P_1=P'_{21}$ and the
$y$-path $P_2=P'_1P'_{22}$ are disjoint and contain all vertices of
$G$.

\begin{subsubcase}
There exists a vertex in $V(H_1)\backslash\{y\}$ with weighted
degree less than $d$.
\end{subsubcase}

Now, we have that $\max\{d^w(v_1),d^w(v_2)\}\geq d$ for every pair
nonadjacent vertices $v_1$ and $v_2$ in $V(H_1)\backslash\{y\}$, and
$\max\{d^w(v_1),d^w(v_2)\}\geq d$ for every pair nonadjacent
vertices $v_1$ and $v_2$ in $V(H_2)$.

Recall that we have assumed that every $(x,z)$-path in $G_2$ has
weight less than $d$. By Theorem 5, there exists a Hamilton
$(x,z)$-path $P'_2$ in $G_2$.

\begin{subsubsubcase}
$H_2$ has at least two vertices.
\end{subsubsubcase}

Let $G'_1$ be the weighted graph such that
$V(G'_1)=V(G_1)\cup\{x'\}$, where $x'\notin V(G)$;
$E(G'_1)=E(G_1)\cup\{xx',x'z\}$; and
$$
w_{G'_1}(e)=\left\{
\begin{array}{ll}
  w_{G_1}(e),   & \mbox{if\ } e\in E(G_1);\\
  0,            & \mbox{otherwise}.
\end{array}
\right.
$$
It is easy to know that $G'_1$ is 2-connected. Besides, we can see
that
$$
\max\{d^w(v_1),d^w(v_2),d^w(v_3)\}\geq d
$$
for every three pairwise nonadjacent vertices $v_1$, $v_2$ and $v_3$
in $V(G'_1)\backslash\{x,y\}$. By the induction hypothesis, $G'_1$
satisfies the conclusion of the theorem.

Suppose that $G'_1$ contains an $(x,y)$-path $P'_1$ of weight at
least $d$. If $P'_1$ contains neither the added edge $xz$ nor the
subpath $xx'z$, then it is an $(x,y)$-path in $G$ with weight at
least $d$. Otherwise, we can use any $(x,z)$-path in $G_2$ to
replace the added edge $xz$ or the subpath $xx'z$, and then obtain
an $(x,y)$-path in $G$ with weight at least $w_{G'_1}(P'_1)\geq d$.

Suppose that $G'_1$ contains an $x$-path $P'_{11}$ and a $y$-path
$P'_{12}$ which are disjoint and have weight sum at least $d$. If
both these two paths contain none of the added edges in
$\{xz,xx',x'z\}$, then they are an $x$-path and a $y$-path which are
disjoint and have weight sum at least $d$. Otherwise, if $P'_{11}$
contains either the added edge $xz$ or one of the subpath in
$\{xx'z,xzx'\}$, then we can use any $(x,z)$-path in $G_2$ to
replace it, and obtain an $x$-path and a $y$-path which are disjoint
and have weight sum at least $d$. Otherwise, $P'_{11}$ or $P'_{12}$
contains the added edge $xx'$ or $zx'$. Then we can use either an
$x$-path in $G_2$ not passing through $z$ or a $z$-path in $G_2$ not
passing through $x$ to replace it, and obtain an $x$-path and a
$y$-path which are disjoint and have weight sum at least $d$.

Suppose that $G'_1$ contains an $x$-path $P'_{11}$ and a $y$-path
$P'_{12}$ which are disjoint and contain all vertices of $G_1$. If
$P'_{11}$ contains the subpath $xx'z$ or $xzx'$, then we can use the
path $P'_2$ to replace the subpath $xx'z$ or $xzx'$, and obtain an
$x$-path and a $y$-path which are disjoint and contain all vertices
of $G$. Otherwise $P'_{11}$ or $P'_{12}$ contains the added edge
$xx'$ or $zx'$, then we can use either $P'_2-z$ or $P'_2-x$ to
replace it, and obtain an $x$-path and a $y$-path which are disjoint
and contain all vertices of $G$.

\begin{subsubsubcase}
$H_2$ has only one vertex.
\end{subsubsubcase}

Let $x'$ be the vertex in $H_2$. Thus we have $x'\in
N(x)\backslash\{y\}$, $d(x')=2$ and $d^w(x')<d$. By a similar proof,
we can obtain that there exists a vertex $y'$ in
$N(y)\backslash\{x\}$ such that $d(y')=2$ and $d^w(y')<d$.

Suppose that $y'=z$. Then $d(z)=2$. If $H_1$ has only one vertex
$y$, then the $x$-path $P_1=xx'z$ and the $y$-path $P_2=y$ are
disjoint and contain all vertices of $G$. Thus, we can assume that
$H_1$ has at least one vertex other than $y$. Let
$G'_1=G[H_1\cup\{x\}]$. If $xy\notin E(G'_1)$, we add the edge $xy$
of weight 0 to $G'_1$, and then get that $G'_1$ is 2-connected and
$\max\{d^w_{G'_1}(v_1),d^w_{G'_1}(v_2)\}\geq d$ for every pair of
nonadjacent vertices $v_1$ and $v_2$ in $V(G'_1)\backslash\{x,y\}$.
By Theorem 5, there exists an $(x,y)$-path $P'_1$ such that either
$w_{G'_1}(P'_1)\geq d$ or $P'_1$ is a Hamilton path in $G'_1$. It is
obvious that $P'_1$ is not the added edge $xy$ itself. If
$w_{G'_1}(P'_1)\geq d$, then $P=P'_1$ is an $(x,y)$-path in $G$ with
weight at least $d$. Otherwise $P'_1$ is a Hamilton path in $G'_1$.
Then the $x$-path $P_1=xx'z$ and the $y$-path $P_2=P'_1-x$ are
disjoint and contain all vertices of $G$.

So, we suppose that $y'\neq z$. Let $y''$ is the vertex adjacent to
$y'$ other than $y$. We can see that each vertex in
$V(G)\backslash\{x,y\}$ other than $x',z,y',y''$ has weighted degree
at least $d$.

Suppose that $y''=z$. Then each vertex in $V(G)\backslash\{x,y\}$
other than $x',z,y'$ has weighted degree at least $d$. Let $G'$ be
the weighted graph obtained by adding the edge $x'y'$ of weight 0 to
$G$. Thus,
$$
\max\{d^w_{G'}(v_1),d^w_{G'}(v_2)\}\geq d
$$
for every pair of nonadjacent vertices $v_1$ and $v_2$ in
$V(G')\backslash\{x,y\}$. By Theorem 5, there exists an $(x,y)$-path
$P'$ in $G'$ such that either $w_{G'}(P')\geq d$ or $P'$ is a
Hamilton path in $G'$. If $P'$ does not contain the added edge
$x'y'$, then $P'$ is an $(x,y)$-path in $G$ with weight at least
$d$, or for every $e\in E(P)$, $P'-e$ is the union of an $x$-path
and a $y$-path which are disjoint and contain all vertices of $G$.
Otherwise, $P'-x'y'$ is an $x$-path and a $y$-path which are
disjoint and either have weight sum at least $d$ or contain all
vertices of $G$.

Suppose that $y''\neq z$. Let $G'_1$ be the weighted graph obtained
by adding two edges $zy'$ and $zy''$ with weight 0 to $G_1$. Thus,
$$
\max\{d^w_{G'_1}(v_1),d^w_{G'_1}(v_2)\}\geq d
$$
for every pair of nonadjacent vertices $v_1$ and $v_2$ in
$V(G'_1)\backslash\{x,y\}$. By Theorem 5, there exists an
$(x,y)$-path $P'_1$ in $G'_1$ such that either $w_{G'_1}(P'_1)\geq
d$ or $P'_1$ is a Hamilton path in $G'_1$.

Suppose that $w_{G'_1}(P'_1)\geq d$. If $P'_1$ contains none of the
added edges in $\{xz,zy',zy''\}$, then the path $P=P'_1$ is an
$(x,y)$-path in $G$ with weight at least $d$. Otherwise we can
delete the added edges in $P'_1$ and then obtain an $x$-path and a
$y$-path which are disjoint and have weight sum at least $d$.

Suppose now that $P'_1$ is a Hamilton path in $G'_1$. If $P'_1$
contains none of the added edges in $\{xz,zy',zy''\}$, or contains
only the added edge $xz$ of them, then the $x$-path $P_1=xx'$ and
the $y$-path $P_2=P'_1-x$ are disjoint and contain all vertices of
$G$. Otherwise, if $P'_1$ contains the added edge $xz$ and one of
the added edges in $\{zy',zy''\}$, then the $x$-path $P_1=xx'z$ and
the $y$-path $P_2=P'_1-\{x,z\}$ are disjoint and contain all
vertices of $G$. So we assume that $P'_1$ does not contains the
added edge $xz$. If $P'_1$ contains one of the added edges $e$ in
$\{zy',zy''\}$, then $P'_1-e$ is the union of an $x$-path and a
$y$-path which are disjoint and contain all vertices in $G'_1$, and
one of these two paths is ending in $z$. Then we can add the edge
$zx'$ to it and obtain an $x$-path and a $y$-path which are disjoint
and contain all vertices of $G$. If $P'_1$ contains both the added
edges $zy'$ and $zy''$, let $P''_1=P'_1\cup\{yy'\}-\{zy',zy''\}$,
then the $x$-path $P_1=xx'z$ and the $y$-path $P_2=P''_1-x$ are
disjoint and contain all vertices of $G$.

The proof is complete.\hfill$\Box$

\section{A related discussion}

If we weaken the condition of Theorem 7 by that of large weighted
degree sum of four pairwise nonadjacent vertices, the result is also
not true.

\begin{false statement}
Let $G$ be a 2-connected weighted graph, $x$ be a vertex of $G$, and
$d$ a be real number. If $d^w(v_1)+d^w(v_2)+d^w(v_3)+d^w(v_4)\geq
4d$ for every four nonadjacent vertices $v_1$, $v_2$, $v_3$ and
$v_4$ in $V(G)\backslash\{x\}$, then $G$ contains either an $x$-path
of weight at least $d$ or a Hamilton $x$-path.
\end{false statement}

This can be shown as follows. Let $H_1$, $H_2$ and $H_3$ be three
complete graphs on at least three vertices with two vertices $x$ and
$y$ in common. Let $G=H_1\cup H_2\cup H_3$ and assign weight 0 to
every edge of $G$. Then $G$ is 2-connected and there exist no four
pairwise nonadjacent vertices in $G$. So $G$ satisfies the condition
of False Statement 3 for any $d>0$, but contains neither an $x$-path
of weight at least $d$ nor a Hamilton $x$-path.

Motivated by Theorems 7 and 8, we pose the following problems.

\begin{problem}
Let $G$ be a 2-connected weighted graph and $d$ a real number. If
$d^w(v_1)+d^w(v_2)+d^w(v_3)+d^w(v_4)\geq 4d$ for every four pairwise
nonadjacent vertices $v_1$, $v_2$, $v_3$ and $v_4$ in $V(G)$, is it
true that $G$ contains a path of weight at least $d$ or a Hamilton
path?
\end{problem}

\begin{problem}
Let $G$ be a 2-connected weighted graph and $d$ a real number. If
$\max\{d^w(v_1)$, $d^w(v_2),d^w(v_3),d^w(v_4)\}\geq d$ for every
four pairwise nonadjacent vertices $v_1$, $v_2$, $v_3$ and $v_4$ in
$V(G)$, is it true that $G$ contains a path of weight at least $d$
or a Hamilton path?
\end{problem}

\end{document}